\documentclass[11pt]{article}

\usepackage{theorem,amssymb,amsmath}

\advance \topmargin by -\headheight
\advance \topmargin by -\headsep
\evensidemargin \oddsidemargin
\author{J.-P. Allouche\footnote{The author was partially 
supported by the ANR project ``FAN'' (Fractal and Numeration).} \\
CNRS, Institut de Math\'ematiques de Jussieu \\
\'Equipe Combinatoire et Optimisation \\
Universit\'e Pierre et Marie Curie, Case 247 \\
4 Place Jussieu \\
F-75252 Paris Cedex 05 France \\
{\tt allouche@math.jussieu.fr}
}

\title{On a formula of T. Rivoal}

\date{ }

\def \proof{\bigbreak\noindent{\it Proof.\ \ }}

\def \endpf{{\ \ $\Box$ \medbreak}}

\newtheorem{theorem}{Theorem}
\newtheorem{lemma}[theorem]{Lemma}
\newtheorem{corollary}[theorem]{Corollary}

\theorembodyfont{\rm}

\newtheorem{remark}[theorem]{Remark}
\newtheorem{example}[theorem]{Example}

\begin{document}

\maketitle

\begin{center}
{\it Dedicated to Professors Zolt\'an Dar\'oczy and Imre K\'atai for their 75th birthday}
\end{center}

\begin{abstract}
In an unpublished 2005 paper T. Rivoal proved the formula
$$
\frac{4}{\pi} = \prod_{k \geq 2} \left(1 + \frac{1}{k+1}\right)^{2\rho(k)\lfloor \log_2(k) - 1 \rfloor}
$$
where $\lfloor x \rfloor$ denotes the (lower) integer part of the real number $x$, and 
$\rho(k)$ is the $4$-periodic sequence defined by $\rho(0) = 1$, $\rho(1) = -1$,
$\rho(2) = \rho(3) = 0$. We show how a lemma in a 1988 paper of J. Shallit and the author
allows us to prove that formula, as well as a family of similar formulas involving occurrences 
of blocks of digits in the base-$B$ expansion of the integer $k$, where $B$ is an integer $\geq 2$.
\end{abstract}

\section{Introduction}

The author shares probably with many number theorists a kind of fascination for infinite products 
or series that look simple, but have explicit and somehow unexpected values, such as
$$
\sum_{n \geq 1} \frac{1}{n^2} = \frac{\pi^2}{6} \ \ \ \mbox{\rm and } \ \ \ 
\prod_{n \geq 1} \left(1 - \frac{1}{4n^2}\right) = \frac{2}{\pi}\cdot
$$
In an unpublished paper \cite{rivoal}, T. Rivoal proved the formula
$$
\frac{4}{\pi} = \prod_{k \geq 2} \left(1 + \frac{1}{k+1}\right)^{2\rho(k)\lfloor \log_2(k) - 1 \rfloor}
$$
where $\lfloor x \rfloor$ denotes the (lower) integer part of the real number $x$, and
where $\rho(k)$ is the $4$-periodic sequence defined by $\rho(0) = 1$, $\rho(1) = -1$,
$\rho(2) = \rho(3) = 0$. Of course this can also be written
$$
\frac{4}{\pi} = \prod_{k \geq 4} \left(1 + \frac{1}{k+1}\right)^{2\rho(k)\lfloor \log_2(k) - 1 \rfloor}.
$$
Grouping terms, this infinite product can also be written
$$
\frac{4}{\pi} = \prod_{k \geq 1} \prod_{0 \leq r \leq 3} 
\left(1 + \frac{1}{4k+r+1}\right)^{2\rho(4k+r)\lfloor \log_2(4k+r)-1 \rfloor}
$$
i.e.,
$$ 
\frac{4}{\pi} = \prod_{k \geq 1} \prod_{0 \leq r \leq 1}
\left(1 + \frac{1}{4k+r+1}\right)^{2\rho(4k+r)\lfloor \log_2(k) + 1 \rfloor}
= \prod_{k \geq 1} \left(\frac{(4k+2)(4k+2)}{(4k+1)(4k+3)}\right)^{2\lfloor \log_2(k) + 1 \rfloor}.
$$
Now, for $k \geq 1$, the quantity $\lfloor \log_2(k) + 1 \rfloor$ is the number of digits in the 
base-$2$ expansion of $k$. Hence, letting $N_{0,2}(k)$ (resp.\ $N_{1,2}(k)$) denote the number 
of occurrences of $0$'s (resp.\ $1$'s) in the binary expansion of the integer $n$, we have 
$\lfloor \log_2(k) + 1 \rfloor = N_{0,2}(k) + N_{1,2}(k)$. Hence Rivoal's relation reads
\begin{equation}\label{eq-rivoal}
\prod_{k \geq 1} \left(\frac{(4k+2)(4k+2)}{(4k+1)(4k+3)}\right)^{2(N_{0,2}(k) + N_{1,2}(k))}
= \frac{4}{\pi}\cdot 
\end{equation}

%
%
%
%

\section{The main result for base $2$}

The purpose of this section is to establish a general relation of which Equation~(\ref{eq-rivoal})
is a particular case. We begin with some definitions. In what follows $B \geq 2$ is an integer
which will be a numeration base for the integers. The set, or {\em alphabet}, ${\cal D}_B$
is defined by ${\cal D}_B := \{0, 1, \cdots, B-1\}$. If $w$ is a {\it word} over ${\cal D}_B$
(i.e., a finite sequence of elements of ${\cal D}_B$), we let $L(w)$ denote its {\em length}: 
if $w = d_1 d_2 \cdots d_k$, then $L(w) = k$ (the usual notation is $|w|$, but $| \ |$ denotes 
the absolute value in a few places in this paper). Also $w^j$ stands for the concatenation of 
$j$ copies of the word $w$.

\medskip

If $w$ is a word over ${\cal D}_B$, we let $N_{w,B}(n)$ denote the number of possibly 
overlapping occurrences of $w$ in the $B$-ary expansion of the integer $n > 0$ if $w$
begins in a $1$ or is of the form $w = 0^j$ for some $j \geq 1$, and the number of possibly 
overlapping occurrences of $w$ in the $B$-ary expansion of the integer $n > 0$ preceded
by an arbitrarily large number of $0$'s if $w$ begins in a $0$, but is not of the form $w = 0^j$ 
for some $j \geq 1$. Finally we define $N_{w,B}(0) = 0$ for any $w$ (which means that $0$ is 
represented by the empty word in base $B$).

\medskip

If $w$ and $B$ are as above, we let $v_B(w)$ denote the ``value'' of $w$ when $w$
is interpreted as the base $B$-expansion (possibly with leading $0$'s) of an
integer.

\begin{example}
To make the above definitions clear we give the following examples:
$N_{11,2}(15) = 3$, $N_{001,2}(4) = 1$ (write $4$ in base $2$ as $0\cdots0100$), while
$N_{0,4}(4) = 2$. Also $v_2(0010) = 2$.
\end{example}

\bigskip

Now we state a general lemma from \cite{as}. A proof is given in \cite{as}
(also see \cite{allshason}, where this lemma is used for proving families
of relations involving the quantities $N_{w,B}(k)$). 

\begin{lemma}[\cite{as}]\label{as}
Fix an integer $B \geq 2$, and let $w$ be a non-empty word over
$\{0, 1, \cdots, B-1\}$. If $f: {\mathbb N} \to {\mathbb C}$ is a function
such that $\sum_{n \geq 1} |f(n)| \log n < \infty$, then
$$
\sum_{n \geq 1} N_{w,B}(n) \left(f(n) - \sum_{0 \leq k \leq B-1} f(Bn+k)\right)
= \sum f(B^{L(w)} n + v_B(w)),
$$
where the last summation is over $n \geq 1$ if $w=0^j$ for some $j \geq 1$,
and over $n \geq 0$ otherwise.
\end{lemma}

\begin{remark}
Note that the relation in Lemma~\ref{as} above does not involve the value $f(0)$.
\end{remark}

The next classical lemma will prove useful (see, e.g., \cite[Section~12-13]{ww}).

\begin{lemma}\label{gamma}
Let $d$ be a positive integer. Let $(a_i)_{1 \leq i \leq d}$ and $(b_j)_{1 \leq j \leq d}$ 
be complex numbers such that no $a_i$ and no $b_j$ belongs to $\{0, -1, -2, \ldots\}$. 
If $a_1 + a_2 + \cdots + a_d = b_1 + b_2 + \cdots + b_d$, then 
$$
\prod_{n \geq 0}\frac{(n+a_1) \cdots (n+a_d)}{(n+b_1) \cdots (n+b_d)}
= \frac{\Gamma(b_1) \cdots \Gamma(b_d)}{\Gamma(a_1) \cdots \Gamma(a_d)}\cdot
$$
\end{lemma}

\begin{theorem}\label{main}
Let $w$ be a word over the alphabet $\{0, 1\}$, and $N_{w,2}$ as defined previously.
Then
\begin{itemize}
\item if $w = 0^j$ for some $j \geq 1$,
$$
\prod_{n \geq 1} \left(\frac{(4n+2)(4n+2)}{(4n+1)(4n+3)}\right)^{2 N_{w,2}(n)}
= \frac{2^{j+2}\Gamma\left(\frac{1}{2^j}\right)}{\Gamma\left(\frac{1}{2^{j+1}}\right)^2};
$$
\item if $w$ is not of the form $0^j$ for some $j \geq 1$,
$$
\prod_{n \geq 1} \left(\frac{(4n+2)(4n+2)}{(4n+1)(4n+3)}\right)^{2 N_{w,2}(n)} =
\frac{\Gamma\left(\frac{v_2(w)}{2^{L(w)}}\right) \Gamma\left(\frac{v_2(w)+1}{2^{L(w)}}\right)}
     {\Gamma\left(\frac{2v_2(w)+1}{2^{L(w)+1}}\right)^2}\cdot
$$
\end{itemize}
\end{theorem}

\proof Define $f$ by $f(0) = 0$ and for all $n \geq 1$
$$
f(n) := \log \left(\frac{(2n+1)^2}{2n(2n+2)}\right).
$$
Then, applying Lemma~\ref{as} with $B = 2$ and $w$ a word over $\{0, 1\}$, yields
$$
\sum_{n \geq 1} N_{w,2}(n) \left(f(n) - f(2n) - f(2n+1)\right)
= \sum f(2^{L(w)} n + v_2(w))
$$
where the last summation is over $n \geq 1$ if $w=0^j$ for some $j \geq 1$,
and over $n \geq 0$ otherwise. Since 
$$
f(n) - f(2n) - f(2n+1) = \log\left(\frac{(4n+2)^4}{(4n+1)^2(4n+3)^2}\right) =
2\log\left(\frac{(4n+2)(4n+2)}{(4n+1)(4n+3)}\right),
$$
we get
$$
\begin{array}{ll}
\displaystyle\sum_{n \geq 1} 2 N_{w,2}(n) \log\left(\frac{(4n+2)(4n+2)}{(4n+1)(4n+3)}\right) =&
\\
& \!\!\!\!\!\!\!\!\!\!
  \!\!\!\!\!\!\!\!\!\!
  \!\!\!\!\!
\displaystyle\sum \log \left(\frac{(2^{L(w)+1} n + 2v_2(w) + 1)(2^{L(w)+1} n + 2v_2(w) + 1)}
{(2^{L(w)+1} n + 2v_2(w))(2^{L(w)+1} n + 2v_2(w) + 2)}\right)
\end{array}
$$
where again the last summation is over $n \geq 1$ if $w=0^j$ for some $j \geq 1$,
and over $n \geq 0$ otherwise. Exponentiating yields
$$
\prod_{n \geq 1} \left(\frac{(4n+2)(4n+2)}{(4n+1)(4n+3)}\right)^{2 N_{w,2}(n)}
= \prod \left(\frac{(2^{L(w)+1} n + 2v_2(w) + 1)(2^{L(w)+1} n + 2v_2(w) + 1)}
{(2^{L(w)+1} n + 2v_2(w))(2^{L(w)+1} n + 2v_2(w) + 2)}\right)
$$
where the last product is over $n \geq 1$ if $w=0^j$ for some $j \geq 1$,
and over $n \geq 0$ otherwise.
Using Lemma~\ref{gamma} (recall that the range of summation for the sum on the 
right of the formula in that lemma is not the same for $w = 0^j$ and for $w \neq 0^j$), 
we then get the statement of the theorem. \endpf

\begin{corollary}
Equation~\ref{eq-rivoal} holds.
\end{corollary}

\proof Applying Theorem~\ref{main} first with $w=0$, then with $w=1$, we obtain
(note that $v_2(0) = 0$, $v_2(1) = 1$, and remember that $\Gamma(1+x) = x\Gamma(x)$)
$$
\prod_{n \geq 1} \left(\frac{(4n+2)(4n+2)}{(4n+1)(4n+3)}\right)^{2 N_{0,2}(n)}
= \frac{8\Gamma\left(\frac{1}{2}\right)}{\Gamma\left(\frac{1}{4}\right)^2}
$$
and
$$
\prod_{n \geq 1} \left(\frac{(4n+2)(4n+2)}{(4n+1)(4n+3)}\right)^{2 N_{1,2}(n)} =
\frac{\Gamma\left(\frac{1}{2}\right)}
     {\Gamma\left(\frac{3}{4}\right)^2}\cdot
$$
Thus
$$
\prod_{n \geq 1} \left(\frac{(4n+2)(4n+2)}{(4n+1)(4n+3)}\right)^{2(N_{0,2}(n) + N_{1,2}(n))}
= \frac{8\Gamma\left(\frac{1}{2}\right)^2}
       {\Gamma\left(\frac{1}{4}\right)^2\Gamma\left(\frac{3}{4}\right)^2}\cdot
$$
But, using Euler's reflection formula $\Gamma(z)\Gamma(1-z)= \pi/\sin(\pi z)$ 
(see, e.g., \cite[Section~12-14]{ww}), we get the classical relations 
$\Gamma(1/2) = \sqrt{\pi}$ and $\Gamma(1/4)\Gamma(3/4) = \pi \sqrt{2}$,
which finally yield 
$$
\prod_{n \geq 1} \left(\frac{(4n+2)(4n+2)}{(4n+1)(4n+3)}\right)^{2(N_{0,2}(n) + N_{1,2}(n))}
= \frac{4}{\pi},
$$
i.e., Equation~(\ref{eq-rivoal}). \endpf

\begin{remark}
We note that the proof of Theorem~\ref{main} gives a companion formula to
Equation~(\ref{eq-rivoal}), namely
\begin{equation}
\prod_{k \geq 1} \left(\frac{(4k+2)(4k+2)}{(4k+1)(4k+3)}\right)^{2(N_{0,2}(k) - N_{1,2}(k))}
= \frac{8 \Gamma\left(\frac{3}{4}\right)^2}{\Gamma\left(\frac{1}{4}\right)^2}
= \frac{16\pi^2}{\Gamma\left(\frac{1}{4}\right)^4},
\end{equation}
but that we were unable to compute the infinite ``alternate'' product 
(see Section~\ref{alternate} for a motivation)
$$
\prod_{k \geq 1} \left(\frac{(4k+2)(4k+2)}{(4k+1)(4k+3)}\right)^{2(-1)^k(N_{0,2}(k) + N_{1,2}(k))}\cdot
$$
\end{remark}

\section{A few words about generalizations to base $B$}

It is actually possible to obtain formulas similar to Rivoal's for bases $B$, 
where $B > 2$. For example Theorem~\ref{main} can be generalized as follows.

\begin{theorem}\label{main2}
Let $w$ be a word over the alphabet $\{0, 1, \ldots, B-1\}$.
Let $(a_i)_{1 \leq i \leq d}$ and $(b_j)_{1 \leq j \leq d}$ be 
nonnegative real numbers.
If $a_1 + a_2 + \cdots + a_d = b_1 + b_2 + \cdots + b_d$, then 

\begin{itemize}

\item if $w = 0^j$ for some $j \geq 1$,

$$
\hskip -1truecm
\prod_{n \geq 1} \left(\prod_{1 \leq i \leq d} \left(\left(\frac{Bn+a_i}{Bn+b_i}\right)
\prod_{0 \leq k \leq B-1} \left(\frac{B^2 n + B k + b_i}{B^2 n + B k + a_i}\right)\right)\right)^{N_{w,B}(n)}
= \prod_{1 \leq i \leq d} \frac{\Gamma\left(1 + \frac{b_i}{B^{j+1}}\right)}
{\Gamma\left(1+ \frac{a_i}{B^{j+1}}\right)};
$$

\item if $w$ is not of the form $0^j$ for some $j \geq 1$,

$$
\hskip -1truecm
\prod_{n \geq 1} \left(\prod_{1 \leq i \leq d} \left(\left(\frac{Bn+a_i}{Bn+b_i}\right)
\prod_{0 \leq k \leq B-1} \left(\frac{B^2 n + B k + b_i}{B^2 n + B k + a_i}\right)\right)\right)^{N_{w,B}(n)}
=
\prod_{1 \leq i \leq d} \frac{\Gamma\left(\frac{v_B(w)}{B^{L(w)}} + \frac{b_i}{B^{L(w)+1}}\right)}
{\Gamma\left(\frac{v_B(w)}{B^{L(w)}} + \frac{a_i}{B^{L(w)+1}}\right)}\cdot
$$

\end{itemize}

\end{theorem}

\proof Apply Lemma~\ref{as} with $f$ defined by $f(0) = 0$ and for all $n \geq 1$
$$
f(n) := \log \prod_{1 \leq i \leq d} \frac{Bn+a_i}{Bn+b_i}\cdot 
$$
\endpf

\begin{remark}
Theorem~\ref{main2} contains Theorem~\ref{main} (take $B=2$, $a_1=a_2=1$, 
$b_1=0$, and $b_2=2$).
\end{remark}

\section{Conclusion}

\subsection{The ``alternate'' Euler constant}\label{alternate}

When he obtained Equation~(1), or more precisely the formula
$$
\frac{4}{\pi} = \prod_{k \geq 2} \left(1 + \frac{1}{k+1}\right)^{2\rho(k)\lfloor \log_2(k) - 1 \rfloor}
$$
Rivoal was inspired by Catalan's and Vacca's identities for the Euler-Mascheroni constant $\gamma$
$$
\gamma = \int_0^1 \frac{\sum_{n \geq 1} x^{2^n}}{x(1+x)}\mbox{\rm d}x
\ \ \ \mbox{\rm and} \ \ \
\gamma = \sum_{k \geq 1} (-1)^k\frac{\lfloor \log_2(k) \rfloor}{k}
$$
(Catalan's identity dates back to 1875, see \cite{catalan}, while Vacca's identity was 
proved in 1925, see \cite{vacca}; for a history of the second formula, see \cite{sondow2}).
An analogy between $\gamma$ and $\log\frac{4}{\pi}$ occurs when writing the above
relations as
$$
\gamma = \sum_{k \geq 1} (-1)^k\frac{\lfloor \log_2(k) \rfloor}{k}
\ \ \ \mbox{\rm and} \ \ \
\log \frac{4}{\pi} = \sum_{k \geq 1} (2\rho(k)\lfloor \log_2(k) - 1 \rfloor)
\log\left(1 + \frac{1}{k+1}\right)\cdot
$$
Another similarity is given by the formulas
$$
\gamma = \sum_{j \geq 2} \frac{(-1)^j}{j}\zeta(j) \ \ \mbox{\rm and} \ \
\log \frac{4}{\pi} = \sum_{j \geq 2} \frac{(-1)^j}{j}\eta(j)
$$ 
where $\eta(j) := (1 - 2^{1-j})\zeta(j)$ (we use the same notation as, e.g., in 
\cite{elsner-prevost} where several formulas of the same kind can be found; also
see \cite{choi-srivastava} and \cite{sondow1}):
these formulas can be obtained by taking $z=1$ and $z=1/2$ in the relation
$$
\log \Gamma(1+z) = - \log(1+z) + z(1-\gamma) + \sum_{n \geq 2} \frac{(-1)^n (\zeta(n) - 1)}{2^n n} 
$$
valid for $|z| < 2$, see \cite[6.1.33, p.\ 256]{abramowitz-stegun}, which gives respectively
$$
\gamma = \sum_{j \geq 2} \frac{(-1)^j}{j}\zeta(j) \ \ \mbox{\rm and} \ \
\gamma = \log \frac{4}{\pi} + 2 \sum_{j \geq 2} \frac{(-1)^j \zeta(j)}{2^j j}\cdot
$$
A more striking analogy between the constants $\gamma$ and $\log \frac{4}{\pi}$ was
noted by Sondow in \cite{sondow1} where it is proved that
$$
\gamma = \sum_{n \geq 1} \left(\frac{1}{n} - \log \frac{n+1}{n}\right)
       = \int\!\!\!\!\int_{[0,1]^2} \frac{1-x}{(1 - xy)(- \log xy)} \mbox{\rm d}x\mbox{\rm d}y
$$
and
$$
\log\frac{4}{\pi} = \sum_{n \geq 1} (-1)^{n-1} \left(\frac{1}{n} - \log \frac{n+1}{n}\right)
       = \int\!\!\!\!\int_{[0,1]^2} \frac{1-x}{(1 + xy)(- \log xy)} \mbox{\rm d}x\mbox{\rm d}y
$$
leading Sondow to call ``alternating Euler constant'' the quantity $\log\frac{4}{\pi}$.
In the same spirit Sondow compares in \cite{sondow2} the following two expressions
$$
\gamma = \frac{1}{2} + \sum_{n \geq 1, \ {\rm even}} \frac{N_{1,2}(n) + N_{0,2}(n) - 1}{n(n+1)(n+2)}
\ \ \ \mbox{\rm and} \ \ \
\log \frac{4}{\pi} = 
\frac{1}{4} + \sum_{n \geq 1, \ {\rm even}} \frac{N_{1,2}(n) - N_{0,2}(n)}{n(n+1)(n+2)}
$$
where the first expression is due to Addison \cite{addison} and the second is a modification of a 
formula in \cite{allshason}.

\bigskip

One way of ``explaining'' the links between $\gamma$ and $\log 4/\pi$ is the introduction of
the ``generalized-Euler-constant function'' by Sondow and Hadjicostas in \cite{sh}, or of a
similar function introduced by Pilehrood and Pilehrood in \cite{pp1}: the function $\gamma(z)$
of \cite{sh} and the function $f_1(z)$ in \cite{pp1} are defined by
$$
\gamma(z) = \sum_{n \geq 1} z^{n-1} \left(\frac{1}{n} - \log \frac{n+1}{n}\right)
\ \ \ \mbox{\rm and} \ \ \ 
f_1(z) = \sum_{n \geq 1} z^n \left(\frac{1}{n} - \log \frac{n+1}{n}\right)
$$
(so that $f_1(z) = z \gamma(z)$). Namely one has
$$
\gamma = \gamma(1) 
\ \ \ \mbox{\rm and} \ \ \
\log\left(\frac{4}{\pi}\right) = \gamma(-1)
$$
(for more on $\gamma(z)$ see \cite{pp2}).

\subsection{Catalan-type formulas}

In his paper \cite{rivoal} Rivoal gives a Catalan-like formula for $4/\pi$
in relation with Equality~\ref{eq-rivoal}, namely 
$$
\sum_{k \geq 2} (2\rho(k)\lfloor \log_2(k) - 1 \rfloor) \log \left(1 + \frac{1}{k+1}\right)
= \log \frac{4}{\pi} 
= \int_0^1 \frac{x-1}{\log x} \frac{\sum_{n \geq 2} x^{2^n}}{x(1+x)(1+x^2)}\mbox{\rm d}x.
$$
Comparing with Catalan's identity
$$
\gamma = \int_0^1 \frac{\sum_{n \geq 1} x^{2^n}}{x(1+x)}\mbox{\rm d}x,
$$
Rivoal suggested (private communication) that similar relations may exist for logarithms 
of the infinite products we studied here. However we do not have general results in that 
direction.

\subsection{Two more remarks}

We would like to make two more remarks about Lemma~\ref{as}.

\begin{itemize}

\item
Which functions can be obtained on the left side of the equality given in that lemma?
In other words given a map $g$ from the integers to the real numbers, we want to know
when it is possible to find a map $f$ such that
$$
g(n) = f(n) - \sum_{0 \leq j \leq B-1} f(Bn+j).
$$
A particular case is addressed in \cite{as}, the case where $f$ is a constant multiple of $g$.
In other words what are the eigenvectors of the operator $f: \to Tf$, where 
$Tf(n) = \sum_{0 \leq j \leq B-1} f(Bn+j)$, and $f$ is supposed to behave ``regularly''?
This looks like a functional equation with means: $\sum_{0 \leq j \leq B-1} f(Bn+j)$ is $B$
times the arithmetic mean of the values of $f$ on $[Bn, Bn+B-1]$. Looking in the literature
for papers with keywords ``mean'' and ``functional equation'', we found several papers,
in particular by Dar\'oczy and coauthors, e.g., \cite{daroczy}, but were not able to find
references really related to our question.

\item
Another question about Lemma~\ref{as} is whether the quantities $N_{w,B}(n)$ can be replaced by
more general sequences. We think that it is possible to introduce generalizations of 
$B$-additive sequences for which a similar lemma holds. We hope to address that question in
the near future, possibly including distribution results (see, e.g., the survey of K\'atai
\cite{katai}).

\end{itemize}

\bigskip

\noindent
{\bf Acknowledgements.} We would like to thank T. Rivoal and J. Shallit for their 
comments on a first version of this paper.


\begin{thebibliography}{99}

\bibitem{abramowitz-stegun} M. Abramowitz, I. Stegun, eds., {\it Handbook of Mathematical 
Functions with Formulas, Graphs, and Mathematical Tables}, National Bureau of Standards, 
Applied Mathematics Series {\bf 55} (1964).

\bibitem{addison} A. W. Addison, A series representation for Euler's constant, {\it Amer. Math. Monthly\, }
{\bf 74} (1967) 823--824.

\bibitem{as} J.-P. Allouche, J. Shallit, Sums of digits and the Hurwitz zeta function,
in Analytic number theory (Tokyo, 1988), Lecture Notes in Math. {\bf 1434}, Springer,
Berlin, 1990, pp. 19--30.

\bibitem{allshason} J.-P. Allouche, J. Shallit, J. Sondow, Summation of series defined by 
counting blocks of digits, {\it J. Number Theory\,} {\bf 123} (2007) 133--143.

\bibitem{catalan} E. Catalan, Sur la constante d'Euler et la fonction de Binet, 
{\it J. Liouville [J. Math. Pures Appl.]} (3) I (1875) 209--241.

\bibitem{choi-srivastava} J. Choi, H. M. Srivastava, Sums associated with the Zeta function,
{\it J. Math. Anal. Appl.} {\bf 206} (1997) 103--120.

\bibitem{daroczy} Z. Dar\'oczy, Mean values and functional equations,
{\it Differ. Equ. Dyn. Syst.} {\bf 17} (2009) 105--113.

\bibitem{elsner-prevost} C. Elsner, M. Pr\'evost, Expansion of Euler's constant in terms of 
Zeta numbers {\it J. Math. Anal. Appl.} {\bf 398} (2013) 508--526.

\bibitem{katai} I. K\'atai, Distributions of arithmetical functions. Some results and problems,
{\it Ann. Univ. Sci. Budapest. Sect. Comput.} {\bf 33} (2010) 239--259. 

\bibitem{pp1} K. H. Pilehrood, T. H. Pilehrood, Arithmetical properties of some series with 
logarithmic coefficients, {\it Math. Z.} {\bf 255} (2007) 117--131 (2007).

\bibitem{pp2} K. H. Pilehrood, T. H. Pilehrood, Vacca-type series for values of the 
generalized-Euler-constant function and its derivative, {\it J. Integer Seq.} {\bf 13} (2010) 
Article 10.7.3.

\bibitem{rivoal} T. Rivoal, Polyn\^omes de type Legendre et approximations de la constante d'Euler,
Unpublished preprint (2005), available at the URL \newline 
{\tt http://www-fourier.ujf-grenoble.fr/$\sim$rivoal/articles/euler.pdf}

\bibitem{sondow1} J. Sondow, Double integrals for Euler's constant and ln($4/\pi$) and an analog 
of Hadjicostas's formula, {\it Amer. Math. Monthly\ } {\bf 112} (2005) 61--65.

\bibitem{sondow2} J. Sondow, New Vacca-type rational series for Euler's constant and its ``alternating'' 
analog $\log 4/\pi$, in {\it Additive Number Theory, Festschrift in Honor of the Sixtieth Birthday of 
M. B. Nathanson}, D. Chudnovsky and G. Chudnovsky, eds., Springer, 2010, pp. 331-340.

\bibitem{sh} J. Sondow, P. Hadjicostas, The generalized-Euler-constant function $\gamma(z)$
and a generalization of Somos's quadratic recurrence constant, {\it J. Math. Anal. Appl.} {\bf 332} 
(2007) 292--314.

\bibitem{vacca} G. Vacca, A new series for the Eulerian constant $\gamma = \cdot 577\ldots$, 
{\it Quart. J. Math} {\bf 41} (1909) 363--364.

\bibitem{ww} E. T. Whittaker, G. N. Watson, {\it A Course of Modern Analysis},
Fourth Edition, reprinted, Cambridge University Press, Cambridge, 1996.


\end{thebibliography}
\end{document}